\documentclass[11pt]{article}

\usepackage{amsthm, amsmath, amssymb,latexsym}

\newtheorem{theorem}{Theorem}[subsection]

\newtheorem{lemma}[theorem]{Lemma}
\newtheorem{corollary}[theorem]{Corollary}
\newtheorem{proposition}[theorem]{Proposition}
\newtheorem{remark}[theorem]{Remark}
\newtheorem{definition}[theorem]{Definition}

\newcommand{\nc}{\newcommand}
\nc{\cH}{{\mathcal H}}
\nc{\cA}{{\mathcal A}}
\nc{\cG}{{\mathcal G}}
\nc{\cC}{{\mathcal C}}
\nc{\cO}{{\mathcal O}}
\nc{\cI}{{\mathcal I}}
\nc{\cB}{{\mathcal B}}
\nc{\cD}{{\mathcal D}}
\nc{\cY}{{\mathcal Y}}
\nc{\cK}{{\mathcal K}} 
\nc{\cX}{{\mathcal X}}
\nc{\cR}{{\mathcal R}}
\nc{\cS}{{\mathcal S}}
\nc{\cE}{{\mathcal E}}
\nc{\cF}{{\mathcal F}}
\nc{\cZ}{{\mathcal Z}}
\nc{\cQ}{{\mathcal Q}}
\nc{\cN}{{\mathcal N}}
\nc{\cP}{{\mathcal P}}
\nc{\cL}{{\mathcal L}}
\nc{\cM}{{\mathcal M}}
\nc{\cT}{{\mathcal T}}
\nc{\cW}{{\mathcal W}}
\nc{\cU}{{\mathcal U}}
\nc{\cJ}{{\mathcal J}}
\nc{\cV}{{\mathcal V}}
\nc{\bH}{{\mathbb H}}
\nc{\bA}{{\mathbb A}}
\nc{\bG}{{\mathbb G}}
\nc{\bC}{{\mathbb C}}
\nc{\bO}{{\mathbb O}}
\nc{\bI}{{\mathbb I}}
\nc{\bB}{{\mathbb B}}
\nc{\bY}{{\mathbb Y}}
\nc{\bK}{{\mathbb K}} 
\nc{\bX}{{\mathbb X}}
\nc{\bS}{{\mathbb S}}
\nc{\bE}{{\mathbb E}}
\nc{\bF}{{\mathbb F}}
\nc{\bZ}{{\mathbb Z}}
\nc{\bQ}{{\mathbb Q}}
\nc{\bN}{{\mathbb N}}
\nc{\bP}{{\mathbb P}}
\nc{\bL}{{\mathbb L}}
\nc{\bM}{{\mathbb M}}
\nc{\bT}{{\mathbb T}}
\nc{\bW}{{\mathbb W}}
\nc{\bU}{{\mathbb U}}
\nc{\bD}{{\mathbb D}}
\nc{\bJ}{{\mathbb J}}
\nc{\bV}{{\mathbb V}}
\nc{\bbZ}{{\mathbb Z}}
\nc{\bR}{{\mathbb R}}
\nc{\co}{{\nabla}}
\nc{\cu}{{\overline{\nabla}}}
\nc{\fr}{{\rightarrow}}
\newcommand{\la}{\longrightarrow}

\begin{document}

\title{On rational maps from a general \\surface in $\mathbb P^3$ 
to surfaces of general type} 
                               
\author{
Lucio Guerra  {\small $^{1}$}
\ and \ Gian Pietro Pirola {\small $^{2}$}
}
\date{}     

\maketitle

\footnotetext[1]{Partially supported by
Finanziamento Ricerca di Base 2006 Univ. Perugia.}

\footnotetext[2]{Partially supported by
1) PRIN 2005 {\em ``Spazi di moduli e teorie di Lie"}; 
2) Indam (GNSAGA);
3) Far 2006 (PV):{\em ``Variet\`{a} algebriche, calcolo
algebrico, grafi orientati e topologici"}.}

\begin{abstract}
{\noindent 
We study dominant rational maps from a general surface in $\mathbb P^{3}$
to  surfaces of general type. We prove restrictions on the target surfaces,
and special properties of these rational maps. 
We show that for small degree
the general surface has no such map.
Moreover a slight improvement of a result of Catanese,
on the number of moduli of a surface of general type,
is also obtained.}

\end{abstract}

\section{Introduction}
 
 Let $X$ be a smooth complex
 projective variety of general 
 type. Let $R(X)$ be the field of rational 
 functions of $X.$ Consider 
 the set  of the {geometric subfields}
 of $R(X)$, that is:
  $$\cF(X)=\{K:\mathbb C \subset K 
 \subset R(X)\}.$$
 An element $K \in \cF(X)$ 
 corresponds 
 to a dominant rational map  $X\dasharrow Y,$
 where $Y$ is a smooth projective variety with  $K\cong R(Y),$ 
 up to birational isomorphisms of $Y.$
 Consider moreover the subset:
 $$\cF_{0}(X)=\{ K \in \cF(X):  [R(X):K] \mbox{ is finite}\}.$$
Elements of this subset correspond to
 generically finite dominant rational maps.
We may then define various 
 {geometric} subsets of $\cF(X),$ such as:
  $$\mathcal{IS}(X)=\{R(Y)\in \cF_{0}(X) : 
  Y\  {\rm is \ of \  general\ type\}},$$  
  $$\cG(X)=\{R(Y)\in \cF_{0}(X) : Y\  
  {\rm is\  not\  rationally \ connected}\}.$$ 
We call $\mathcal{IS}(X)$ the Iitaka-Severi set of 
$X,$ we denote by $s(X)$ 
the cardinality of $\mathcal{IS}(X)$ and by $g(X)$
 the cardinality of $\cG(X).$
 
The recent solution of the Iitaka-Severi conjecture
 \cite{Tsuji,Haconk, takayama} 
 gives that $\mathcal{IS}(X)$ 
is a finite set. In general 
$\cG(X)$  is not finite (for 
instance if $X$ dominates an abelian 
variety).
The problem remains how to compute or 
at least to estimate the number s(X). We call
this the refined Iitaka-Severi problem.

When  $X$ is a
 curve of genus $g\geq 2$ 
 effective bounds on $s(X)$ 
 in terms of $g$ are known \cite{Kani,Tanabe}.
 In higher dimensions not much is known,
but see \cite{Guerra,Heier, NaranjoPirola} where  upper
bounds are given under some geometric restrictions.

If  $X$ is a curve,  general in moduli,
then $g(X)=s(X)=1.$ 
This may be proved by
 counting moduli of maps by means of 
the Hurwitz formula.  
The same in fact
 holds for the general smooth plane 
curve of degree bigger than $3,$  
or for a general
hyperplane section  of a regular surface \cite{ciro, Severi}.
The proof of these facts may be based on
a Hodge-Lefschetz theoretical
 argument (using monodromy, see \cite{Voisin}, \S 3.2.3),
 which implies 
that the Jacobian of $H$ is simple.
So we have two basic methods: 
a moduli count and a Hodge theoretical argument.
  
In higher dimensions,
we believe the following could be true:
\medskip
 
\noindent{\bf Conjecture.}
{\em Let $X$ be a  very general hypersurface of 
$\mathbb P^{n}$ of degree $d>n+1.$
 Then $g(X)=s(X)=1.$}
\medskip

The case of curves and the results of Amerik
 \cite{Amerik} give  evidence to the  conjecture.
 In this paper we can prove the following  (see \ref{main}):  
 \medskip
 
\noindent{\bf Theorem.}   
{\em  If  $X$ is a general surface of $\mathbb P^{3}$
 of degree   $5\leq d\leq 11$ then $s(X)=1.$}
\medskip

The  proof uses both methods described
for the case of curves.
Using the Hodge theoretic argument,
we obtain
restrictions for the target surfaces  (see \ref{hodge}):
\medskip

\noindent{\bf Proposition.} 
{\em Let  $X$ be a general surface of 
 $\mathbb P^{3}$ of degree $d\geq 5.$
Let $Y\!$ be a minimal surface of general type, and
assume that $f: X\dasharrow Y$ is a dominant rational map.
 If $f$ is not birational then $Y$ is simply connected 
of geometric genus $p_g(Y)=0$.} 
\medskip

Here we mention that
simply connected surfaces
of geometric genus $p_g(Y)=0$
are known to exist \cite{Barlow, LeePark}, and moreover 
they are homeomorphic to rational surfaces,
as follows from Freedman's theorem \cite{freedman}. 
However the moduli space of these surfaces 
is still largely unknown.

Then we approach the moduli of rational maps.
First we consider the moduli of target surfaces.
There is a well known result of Catanese 
\cite {Catanese,Catanese2}
on the moduli of surfaces of general type, 
for which we propose a
new approach based on the stability
\cite{bogomolov, enoki, reid, sugiyama},
which produces a slight improvement  (see \ref{mi1}):
\medskip

\noindent{\bf Theorem.}  
{\em Let $Y$ be a minimal surface of general type,  $M(Y)$ be the number
of moduli of $Y.$ The following estimate holds:  
$M(Y)\leq 11\chi(\cO_{Y})+K_{Y}^2.$}
\medskip 

 Then we study the moduli of maps
 in terms of their ramification.
Roughly speaking we associate
 to a rational map $f:X\dasharrow Y$
the piece of the ramification 
divisor that is seen on $X$ (not in the exceptional
divisor of the resolution of $f$).  This is a complete
 intersection curve $D$ on the surface $X \subset \bP^{3}.$

As an outcome of the vanishing $p_{g}(Y)=q(Y)=0$
we have  that a property of Cayley-Bacharach type
is enjoyed by the fibers of the rational map,
and this in turn implies some estimate
for the degree (\ref{cortini})  
and the ramification (\ref{stu2}) of the map.

 A lower bound for
the moduli of the ramification divisor $D \subset X$
follows from an argument (\ref{mori}) which combines
the rigidity theorem for rational maps
and the bend and break lemma of Mori theory.
 An upper bound for
the moduli of the curve $D \subset \bP^{3}$
is obtained in terms of 
the degree of this complete intersection.
Finally all the constraints 
force the inequality $d\geq 12.$ 

\vskip 5pt
 
The problem we studied
was a good field for the interplay of
 methods that come from different areas:
 projective and birational geometry, moduli
and stability theory. 
It is possible that a different approach is necessary in order to 
settle the above conjecture. One possibility is to try
a degeneration argument. The results 
of this paper will be useful 
(in the surface case), 
allowing for instance to 
consider a restricted class of 
target varieties.

\section {Surfaces of general type}

In this section we study the number of moduli 
of a curve which deforms in a surface of general type,
and in the last subsection we study
the number of moduli of a surface of general type.

\subsection{Rigidity}
We need the following
 rigidity theorem for rational maps:
 
\begin{theorem} \label{rigidity}
If $X$ and $Y$ are varieties of general type, of the same dimension, 
a dominant rational map
$f: X \dasharrow Y$  admits no non-constant deformation.
\end{theorem}
This was first proved by Kobayashi-Ochiai \cite{KO},
and also follows from the more general statement
known as the Iitaka-Severi conjecture,
nowdays a theorem in virtue of the
recent work of Tsuji \cite{Tsuji},
Hacon-McKernan \cite{Haconk},
and Takayama \cite{takayama},
and the original approach of 
Maehara \cite{maehara}.
An updated account will be presented in 
a forthcoming paper \cite{lggpp}.

\subsection{Bend and break}

We prove a lemma which combines the rigidity theorem and
the basic idea of Mori theory.

\begin{lemma} 
Let $S$ and $B$ be
 smooth connected projective
 surfaces.  Let  $C$ be
 a smooth connected projective curve. Let $F:C\times B \dasharrow S$ 
be a rational map. 
Assume that the family   $F(C_b)$ 
is a  two-dimensional family
of curves on $S.$
Then $S$  is not of general type. \label{mori}
\end{lemma}
\begin{proof}
Assume by contradiction that $S$ is of general type.
We first remark that $F$ is dominant. 
Take the general point  $s\in S$ 
and the general point
 $(t,b)\in F^{-1}(s).$
We remark that  $F^{-1}(s)$ is a curve on $C\times B.$  
Fix a general point $x\in C$ 
and consider the rational map $F_x: B\dasharrow S:$
$$F_x(b)=F(x,b).$$  

We now show that $F_x$ cannot be dominant.
Otherwise  the  rigidity theorem
 \ref{rigidity}
gives $F_t=F_x$ for $t$ belonging to a Zariski open subset 
of  $B.$ This implies then
$F(t,b)=F(x,b),$ 
that is the curve $C_b$
 is contracted to a point. This gives a contradiction.
 
We have then that $F_x(B)$ is a curve 
and hence there is a
curve $D$ on $B$ such that 
$F(x,y)=s$ for all $y\in D.$
The family of curves obtained
by restriction of $F$
defines a map
$$G:C\times D \dasharrow  S$$ 
such that $G(x,y)=s$ for $y\in D.$ 
Using Mori's trick (see \cite{Mori}, thm. 5 and  6) it follows that there is 
a rational curve $R\subset S$
 passing through $s.$ Therefore 
 $S$ cannot be of general type.
\end{proof}

\subsection{Modular dimension} \label{moddim}

A family  of curves parametrized by a nonsingular
variety $U$  is a 
surjective proper morphism $q: \cX \rightarrow U$,  
with  1-dimensional fibers  $X_{t}=q^{-1}(t)$.
We assume that $\cX$ is reduced, in order to avoid
multiple components in the general curve, we also assume that all
components of $\cX$ dominate $U.$ 
A smooth family of curves will be a family of curves for which $q$ is a smooth morphism. 
Let ${\mathfrak M}_{g}$ be the moduli space of 
smooth  connected curves of genus $g$. For any smooth family of 
connected curves of genus $g,$ there is 
a modular map
 $$\mu: U \to {\mathfrak M}_{g}$$
defined by $\mu(t) = [X_{t}]$ (see e.g. \cite{GIT}, ch. 5).
The dimension of the image of this map
is the number of moduli of curves in the family.
For an arbitrary family of curves we define
the number of moduli as follows.

Assume first that $\cX$ is irreducible. There is a  dominant map
$k: W \to U$, 
which is a finite morphism of $W$ onto a Zariski open subset of $U$,
and there is a smooth connected family of curves $\cC \to W$,
together
with a morphism of families 
$$\begin{array}{ccc}
\cC & \overset{j}{\rightarrow} & \cX \\
\downarrow && \downarrow \\
W & \underset{k}{\rightarrow} & U
\end{array}$$
 such that for $t \in U$ the induced morphism
$$\coprod_{k(z)=t} C_z \to X_{t}$$ 
is the normalization map of $X_{t}$.
This family $\cC \to W$
defines a modular map  
$$\mu : W \to \mathfrak M_g.$$ 

\begin{definition} \label{modulardimension}
Let $\cX \rightarrow U$ be a family of curves.
If $\cX$ is irreducible, 
the dimension of the image of the modular map $W \to \mathfrak M_g $ is
by definition the modular dimension of the family:
$$M(\cX/U) = \dim \mu(W).$$
In general, if $\cX = \bigcup \cX_{i}$
is the irreducible  decomposition,
then by definition the modular dimension of the family is:
$$M(\cX/U)=\max_i M(\cX_i/U).$$
\end{definition}

A family of curves in a variety $Y$
is a family of curves $\cX \rightarrow U$ such that
$\cX \subset U \times Y$. 
Over a Zariski open subset $U' \subset U$
the family is flat,  there is the
natural map  $U' \to \cH(Y)$ 
to the Hilbert scheme of $Y$,
sending $t \mapsto X_{t}$, and the dimension of
the image of this map is the dimension of the family, we call it $f$.
We remark that in general  $M(\cX/U)\leq f.$
We can rewrite \ref{mori}:

\begin{proposition} \label{curve}
Let $\cX/U$ be a
 family of curves on a surface of general
type, let $f$ be the dimension of the family.
Then: $$f-1\leq M(\cX/U) \leq f.$$ 
\end{proposition}
\begin{proof}
The fibers of the modular map $ \mu: W\to \mathfrak{M}_g$ 
define families of curves with constant moduli.
 In  a surface of general type
by \ref {mori} they have dimension $\leq 1.$
\end{proof}

  \subsection{Stability}

Let $M$ be a smooth  projective complex variety of dimension $n.$
Let  $A$ be a line bundle on $M.$  We  would like to  recall the notion 
of Mumford-Takemoto semistability of a 
 vector bundle with respect to $A.$  Usually $A$ is  ample or at least 
 nef and big. We will often abuse notation and identify a line bundle 
$A$ with its first Chern class $c_1(A).$
 Let  $E$ be a vector bundle of rank $r$ on $M.$ We say that $E$
 is semistable with respect to $A$ if for any  injective sheaf map 
 $$\phi  :F \la E,$$ where $F$ is a coherent  sheaf of rank $s,$ then:
\begin{equation*} \label{mutare}
 \frac{ c_1(F)\cdot A^{n-1} }{s} \leq \frac{c_1(E)\cdot A^{n-1} }{r} .
 \end{equation*}
  
Let $T_M$  and $\Omega^1_M$ be   the tangent 
and the cotangent bundle of $M.$
From the work of Yau  \cite{ yau,yau2} it follows
that if the canonical bundle is ample then
 $T_M$   is $K_M$-semistable.
  The following  general version
  was proved by Enoki (see  \cite{sugiyama} and \cite{enoki}):
 
 \begin{theorem} \label{eno11}
 Let $M$ be a canonical projective variety (that is, $M$ has
  only canonical singularities and ample canonical divisor) 
  of dimension $n$.   Let $\mu\colon N\to M$ be a smooth resolution. 
  Then the tangent bundle $T_N$ is $\mu^{\ast}K_M$-semistable.
 \end{theorem}

Now let $Y$ be a minimal surface of general type and $\Omega^1_Y$ 
be its cotangent bundle. The canonical model of $Y$ has canonical singularities.
Then  \ref{eno11}  implies that
 every $ \Omega^1_Y(mK_Y)$  is semistable with respect to
$K_{Y}.$ We obtain the following: 

 \begin{corollary}\label{enoki}
Let $Y$ be a minimal surface of general type. If  
$L$ is a line bundle on $Y$ and 
 $L\subset \Omega^1_Y(mK_Y)$ is a sheaf inclusion, then:
$2 K_{Y}L\leq (2m+1)K^2_{Y}.$
\end{corollary}

\begin{remark}  \em
Bogomolov (see for instance \cite{bogomolov} and \cite{reid})   introduced 
a slightly weaker notion of stability, called the $T$-(semi)stability 
in the book of Kobayashi \cite{Kobayashi}, p. 184. 
Bogomolov was able to prove that for a minimal surface of general type
$\Omega^1_Y$ is $T$-semistable. This is
also a straightforward consequence of \ref{enoki}.
It is  likely possible that  \ref{enoki}  follows  from Yau's work,
or from the theory of  Bogomolov,
but presently the authors do not have a reference for that.
\end{remark}

\subsection{Surface moduli estimate}

In this section we prove a bound for the number of moduli $M(Y)$ of a 
minimal surface $Y$ of general type, 
slightly improving a result of Catanese.
The proof relies on   \ref{enoki}.

\begin{theorem} \label{mi1}
Let $Y$ be a minimal surface of general type, and let $M(Y)$ be the number
of moduli of $Y$. The following estimates hold:  
\begin{enumerate}
\item
$ M(Y)\leq 11\chi(\cO_{Y})+K_{Y}^2;$
\item
if $K^{2}_{Y} =1$ then  $ M(Y)\leq 10\chi(\cO_{Y})+1.$
\end{enumerate}
 \end{theorem}

Let $T_Y$ be the tangent bundle of $Y$. 
We have $h^2(T_Y)=
h^0(\Omega^1_Y(K_Y))$,
 $h^0(T_Y)=0$ and   $h^1(T_Y)\geq M(Y)$, by deformation theory.  
From Riemann-Roch we have  
$\chi (T_Y)=2K_Y^2-10\chi(\cO_Y)$, hence:
$$M(Y)\leq h^0(\Omega^1_Y(K_Y)) - \chi (T_Y) = 
h^0(\Omega^1_Y(K_Y)) +10\chi(\cO_Y)-
2K_{Y}^2.$$

 Therefore \ref{mi1} is an immediate consequence
 of the following: 
 
 \begin{proposition}
 In the present setting we have:
 \begin{enumerate}
 \item
 $h^0(\Omega^1_Y(K_Y))\leq  \chi(\cO_{Y})+ 3K^2_{Y};$
 \item
 if $K^{2}_{Y} =1$ then $h^0(\Omega^1_Y(K_Y))\leq  3.$
\end{enumerate}
 \end{proposition}
 
 \begin{proof}
Assume that there is a line bundle
$L\subset \Omega^1_Y(K_Y)$ with 
$h^0(L)=\dim H^0(Y,L)>0.$ Otherwise $h^0(\Omega^1_Y(K_Y))=0$
and the statement is trivially true since  $\chi(\cO_{Y})+ 3K^2_{Y}\geq 4 .$
After saturation of $L,$ we can define 
an exact sequence:
\begin{equation*} \label{satura}
0\la L\la  \Omega^1_Y(K_Y) \la M\otimes I_{\theta}\la 0,
\end{equation*}
where $M$ is a line bundle and $I_{\theta}$ 
 is the ideal of a zero
dimensional scheme $\theta.$ 
We obtain
$$
 h^0(\Omega^1_Y(K_Y))\leq h^0(L)+
h^0(M\otimes I_{\theta})\leq  h^0(L)+h^0(M).
$$
Note that  $L+M=\det(\Omega^1_Y(K_Y))=3K_Y.$
\bigskip

(a) First assume $K^{2}_{Y} \geq 2$.
 We have two subcases:
 \medskip
 
\noindent
(1) $h^0(M)\neq 0$. In this case 
using the multiplication  
$$\mu: H^0(L) \otimes H^0(M)\la H^0(3K_Y)
$$
we obtain by Hopf's lemma 
$$
h^0(L)+h^0(M)-1\leq \dim {\rm Im}(\mu)\leq h^0(3K_{Y})=
\chi(\cO_{Y})+ 3K^2_{Y}.
$$
We know from base point freeness 
and 1-connectedness (\cite{BHPV}, Ch.7 \S\S 5,6)
that $|3K_{Y}|$ contains some
smooth irreducible curve $D$, and therefore the strict inequality
$h^0(L)+h^0(M)-1 <  h^0(3K_{Y})$ holds. In conclusion:
$$
h^0(\Omega^1_Y(K_Y))\leq  \chi(\cO_{Y})+ 3K^2_{Y}.
$$ 

\noindent
(2) $h^0(L)= h^0(\Omega^1_Y(K_Y))$. 
Take a smooth irreducible curve $D$ in $|3K_{Y}|$, as before. We have 
$h^0(L)\leq h^0(L_{D}).$
In fact $L(-D)$ is contained in $T_{Y}(-K_{Y})$, which has no sections.
We can apply Clifford's  theorem and \ref{enoki} to get:
$$ 2(h^0(L_{D})-1)\leq DL=3K_{Y}L\leq \frac{9}{2}K_Y^2,
$$
so finally
$$ h^0(\Omega^1_Y(K_Y))\leq  h^0(L_{D}) \leq 1+\frac{9}{4} K_Y^2\leq 
\chi(\cO_{Y})+ 3K^2_{Y}.
$$
\smallskip

(b) Now consider the case $K^{2}_{Y} =1$.
Let us prove that:
if  $L\subset \Omega^1_Y(K_Y)$ then 
$h^0(L)\leq 1.$
From  \ref{enoki} we have
 $2K_{Y}L\leq 3$ that is
$$K_{Y}L \leq 1.$$
Assume by contradiction $h^0(L)\geq 2$. 
Write $|L|=F+|H|$
where $ F$ is the fixed part of the system 
and $H$ is the free part.
We have $1\leq K_YH\leq K_YF+K_YH=K_{Y}L=1$
 that is $K_YH=1.$ 
 It follows that $H^2$  is  odd (since $H^{2}-K_YH$ is even by Riemann-Roch) 
 and $  \geq 0.$ 
The Hodge index theorem
gives $H^2=1$ and $H\equiv K_{Y}$ numerically. 
It follows that  $h^0(2H)\geq 3.$ But now 
Ramanujam vanishing gives 
$h^1(2H)=h^1(2K_Y)=h^1(-K_Y)=0$
 and $h^2(2H)=0$. That would imply 
$h^0(2H)=h^0(2K_Y)=1+K_Y^2=2$, which gives a contradiction.

Now consider the determinant map 
$c: \wedge {^2} H^0(\Omega^1_Y(K_Y)) \la H^0(3K_Y).$
From the assertion above, the kernel of $c$ does 
not contain any decomposable
non trivial element. Otherwise, if $s_{1}\wedge s_{2} =0$ then
the two sections define a rank 1 subsheaf $L$ of $\Omega^1_Y(K_Y)$
with $h^{0}(L) \geq 2.$
Since  $h^0(3K_Y)=4$ it follows
 that  $h^0(\Omega^1_Y(K_Y))\leq 3.$ 
\end{proof}

\begin{corollary} \label{mo2}
Let $Y$ be a simply connected minimal surface 
of general type with $p_{g}(Y)=0$. We have
\begin{enumerate}
\item
  $M(Y) \leq K^2_{Y}+11\leq 19;$
\item
 if  $K_Y^2=1$ then 
$M(Y)\leq 11.$ 
\end{enumerate}
\end {corollary}
\begin{proof} 
Under the present hypotheses, we have 
$\chi(\mathcal O_Y)=1$. Then  $K_{Y}^2\leq 9$
by the Miyaoka-Bogomolov inequality. 
Moreover by Yau's theorem
if $K_{Y}^2=9$ then $Y$ is not simply connected (and is rigid).
Then the result follows from \ref{mi1}.
 \end{proof}

  \begin{remark} \em
The following  estimate of the number of 
moduli of minimal surfaces 
of general type was given by  Fabrizio 
Catanese (\cite {Catanese} thm.  B, and \cite{Catanese2} thm. 20.6):
 $$  M(Y) \leq 10\chi(\mathcal O_Y)+3K^2_{Y}+18.$$
 The estimate in \ref{mi1} is a slight improvement,
 as is easily seen using the Noether inequality.

\end{remark}

\begin{remark} \em   \label{mi}
 For a surface $X$ of degree $d\geq 5$  
 of $\mathbb P^{3}$ 
 $$M(X)=M(d)=\binom{d+3}{3}-16=\frac{(d+1)(d+2)(d+3)}{6}-16.$$
 \end{remark}

\section{Surfaces of  projective space}

We study rational maps from a surface in $\mathbb P^{3}$
to a surface of general type. Under certain special assumptions,
we obtain some estimate for the degree of the map and 
some control of the ramification. In the last subsection
we prove that for the general surface in $\mathbb P^{3}$ 
the assumptions are
indeed verified.

\subsection{Hurwitz formula} \label{hurwitz}
In this section $X$ will be a smooth surface 
of $\bP^3$ of degree $d>4,$
with Picard group generated by the 
hyperplane section,
$Y$ will be a minimal surface of general type with
$p_g(Y)=0$, and $f:X\dasharrow Y$ will be a 
dominant rational map of degree $m.$
Consider the diagram of maps
\begin{equation} \label{diagram}
\begin{array}{ccc}
&Z& \\
&^{\scriptsize\phi}\!\! \swarrow 
\hspace{10pt} \searrow \!\! ^{\scriptsize  h} \\
X \hspace{-10pt} & \underset{f}{\dasharrow} & \hspace{-10pt} {Y}
\end{array}
\end{equation}
where 
$\phi$ is the blowing 
up which resolves the singularity of 
$f,$ and $h$ is the morphism
which extends  $f,$ so that 
(as rational maps) $ h\circ \phi^{-1} = f .$

Let $E $ be the ramification divisor of
$\phi: Z\to X.$  Every connected component of the support
of $E$ is a connected tree of rational curves.
Let $H$ be the hyperplane divisor of $X.$
Set $L=\phi^{\ast}H.$
Let $K_{X},$ $K_{Z}$ and $K_{Y}$
be the canonical divisors of
$X,$ $Z$ and $Y.$ We have $$K_{X}=(d-4)H.$$
Since the N\'eron-Severi group of $ X$
is generated by the hyperplane $H$
we have that the N\'eron-Severi group
(=Picard group) of $Z$ is
generated by $L$ and the   irreducible components
$E_{i}$ of the support of $E.$
Let $R\subset Z$ be the ramification divisor of $h.$
The Hurwitz formulae give (modulo linear
equivalence):

\begin{equation} \label{he}
K_{Z}=h^{\ast}(K_{Y})+R=
\phi^{\ast}(K_{X})+E=(d-4)L +E
\end{equation}

\noindent Write:

\noindent \begin{equation} \label{can}
h^{\ast}(K_{Y})=rL-W,
\end{equation}
\begin{equation} R= sL+ W +E, \label{ram}
\end{equation}
where: $$W=\sum a_{i}E_{i}.$$
The coefficients
$a_{i}$ and $r, s$ are integers, with $r \geq 0, s \geq 0$ and
$$r+s=d-4.$$
We prove the following:

\begin{lemma} \label{pos}
The divisor $ W=\sum a_{i}E_{i}$ is effective,
that is  $a_{i}\geq 0$ for all $i.$
\end{lemma}

\begin{proof}
Write
$W=A-B$ where $A$ and $B$ are effective 
divisors supported on $E'$ with disjoint 
irreducible components, in particular :
$$A\cdot B\geq 0.$$
Now  $h^{\ast}(K_{Y})=rL-A+B$ is a nef divisor
since $Y$ is minimal surface of general type.
Then, since  $L\cdot B=0,$ we get:
$$ 0\leq B\cdot h^{\ast}(K_{Y})= 
-B\cdot A+ B^{2}\leq B^{2} .$$
 This implies
$B=0$ since  $B$ is contracted by 
$\phi.$   \end{proof}

\begin{lemma}\label{stu}
Using the notation established above, 
we have: 
$$ 
mK_{Y}^{2}\leq r^{2}d.
$$ 
Moreover  $ mK_{Y}^{2}= r^{2}d$ 
holds if and only if $W=0.$ 
In particular we have $r>0.$

\end{lemma}
 \begin{proof}
Since $L^2=H^{2}=d$  we obtain:
 $$ mK_{Y}^{2}=(h^{\ast}K_{Y})^{2}=
 r^{2}H^{2}+W^{2}=r^{2}d+W^{2}\leq 
 r^{2}d.$$
Since $K_{Y}^{2}>0$ we get $r>0.$ 
\end{proof}
 
 \begin{remark} \em
To show that $r>0$ it is enough to assume 
that $Y$ is simply connected with $p_{g}(Y)=0$
and  $K_Y$  nef. 
This implies  clearly $r\geq 0.$ 
If  we assume by contradiction 
$K_Y^2=0$ and $r=0,$ we would
obtain
 $W=0,$ that is $ f^{\ast}(K_{Y})=0.$  But this would give
$f_{\ast}f^{\ast}(K_{Y})=mK_{Y}=0.$ 
Since $Y$ is  simply 
connected  it would follow then $K_{Y}=0,$
and hence  
$p_{g}(Y)=1.$ This is a contradiction.
\end{remark}

\subsection{Cayley-Bacharach condition}
 
 We recall the classical notion of the Cayley-Bacharach
condition (see \cite{griffithsha}).
A set of distinct points 
 $T=\{p_1,\dots, p_{m-1},p_m\}\subset \bP^3$
is in  Cayley-Bacharach position
  with respect to  $\cO(d)$
  if any surface of degree $d$
passing through any subset of $T$ of cardinality $m-1,$
 must contain also the remaining point.
There are many results 
on points in Cayley-Bacharach position
(see \cite{Cheltsov}).
We will use only the following  
elementary lemma:

\begin{lemma}\label{cb}
Assume  $T=\{p_1,\dots, p_{m-1},p_m\}$ in $\mathbb P^3$ is in 
 Cayley-Bacharach position with
  respect to  $\cO(n),$
$n>0$. Then:
\begin{enumerate}
\item $m\geq n+2.$
\item If $n+2\leq m\leq 3n+1$ then  $T$ is 
contained in a plane.
\item If $n+2\leq m\leq 2n+1 $ then $T$
 is contained in a line.
\end{enumerate}
\end{lemma}

\subsection {Trace of holomorphic forms }
Let $X$ be a
 smooth surface of 
 $\mathbb P^{3}$ of degree $d\geq 5$. Assume that
 $ f: X \dasharrow Y$ is a generically finite
 dominant rational map of 
 degree $m=\deg f$.
 We will use the method of  \cite{lopezpirola}.
 Any rational correspondence between $X$ and
  $Y,$ $\Gamma\subset X\times Y$,
 defines a trace map 
 $ tr(\Gamma): H^{2,0}(X)\to H^{2,0}(Y)$
  defined by the composition
 of pull-back and push-down
  $ \pi{_{Y}}\!_{\ast} \pi^{\ast}_X$, where
   $\pi_X: \Gamma\to X$
 and $ \pi_Y: \Gamma\to Y$ 
 are induced by the projections.
When $f:X\dasharrow Y$  is generically finite,
the trace  of $f$
$$ tr(f):  H^{2,0}(X)\to H^{2,0}(Y)$$
 is associated
 to the graph of $f.$
Let  $y$ be a  general point of $Y,$ and assume that
$f$ is \'{e}tale at $y.$ 
Set  $T=\{ p_1,\dots,p_m\}= f^{-1}(y).$
Taking a local coordinate $z$ around
 the point $y$ we define then by 
 pullback coordinates
around any point $p_i\in T=f^{-1}(y).$
 Now if $\omega \in H^{2,0}(X)$ 
 using the parameters defined above 
 as local identification
we get the local trace formula:

\begin{equation*} \label{tracci2}
tr(f)(\omega)_y=\sum_{p_i\in T}\omega_{p_i}.
\end{equation*}

Assuming that $tr(f)(\omega)=0$ and 
that $\omega$ vanishes in $m-1$ points of 
$T,$ it follows that $\omega$ 
must vanish in the remaing one. 
If $tr(f)=0$, and this certainly happens  when $p_{g}(Y)=
 \dim H^{2,0}(Y)=0,$ then
$T $  is in  Cayley-Bacharach position with
  respect to $\cO(d-4).$ In particular we have then:
  
  \begin{proposition}
Let $X$ be a smooth surface in $\bP^{3}$ of degree $d,$  
 and $Y$ be a 
  smooth surface with $p_g(Y)=0.$
 Let $ f: X \dasharrow Y$ be a generically 
  finite  rational  map of degree $m=\deg f.$
  Then the points of the general fiber of $f$ 
  are in Cayley-Bacharach position
with respect to $\cO(d-4).$
  \end{proposition}

\subsection{Degree of maps}
We now prove the following:

 \begin{proposition}\label{cortini}
 Let $X$ be a smooth surface in $\bP^{3}$ of degree $d,$  
 and $Y$ be a non-rational
  smooth surface with $p_g(Y)=q(Y)=0.$
 Let $ f: X \dasharrow Y$ be a generically 
  finite  rational  map of degree $m=\deg f.$
  We have 
  \begin {enumerate}
 \item    $m\geq d-1.$  
 \item  
If  $d>5$ and $X$ does not 
contain rational curves  then $m\geq d.$ 
\end{enumerate}
\end{proposition}
\begin{proof}
\ 1. \  From \ref{cb} it follows
that $m>d-3$.  Assume  by contradiction
  $m= d-2.$ By  \ref{cb} 
 we still have that the points of a general
  fiber of $f$ are on a line.
Let  $S^{k}(X)$ be the $k$-symmetric
product of $X$  and
define the rational map
$Y \dasharrow S^{d-2}(X),$ $y \to f^{-1}(y).$ Taking 
the two residual
points of  the line which contains $f^{-1}(y)$  
  we define a rational map 
$k:Y \dasharrow S^{2}(X).$ The map $k$ is
birational onto its image. This follows 
since generically the two points 
define the line and then 
the fiber of $f.$  
  The main point of \cite{lopezpirola}  is that
   the image of $k$ cannot define
    a correspondence
between $X$ and $Y,$ since otherwise this
 also should be a trace null correspondence.
The analysis of \cite{lopezpirola} 
proves then that $Y$ is 
 birationally isomorphic either to the
product of two curves of $X$ or to the 
$2$-symmetric product of a curve 
in  $X$ or else to a rational ruled 
  surface over a curve of $X.$
Since $Y$ is regular (it is dominated
 by the regular variety $X$)
it would follow that $k(Y)$ is covered 
by rational curves:  $Y$  is rational, and we 
 obtain a contradiction.   
 \vskip 5pt
 
2. \ Assume now by contradiction $m=d-1$ 
and $d>5.$ We will show that $X$ contains a
rational curve. Since $d>5$  the general
 fiber of $f$  is contained in a line (see \ref{cb}).
 Arguing as before we get  a rational map
$Y\dasharrow S^{d-1}(X)$ and by taking the
 residue point on the line we get a map
$g:Y \dasharrow X.$ This map  cannot be 
dominant ($p_{g}(X)>0$).
Then either $g(Y)$ is a point or a
 curve. If $g(Y)=p$ is a point
 the general fibers of $f$ are the 
 general fibers of the 
projection $\pi_{p}:X \dasharrow \mathbb P^{2}$ from $p,$ then $f=\pi_{p}$ 
(as rational maps) and hence $Y$ is birational to $\mathbb P^{2}.$
This gives a contradiction.
It follows that $g(Y)$ is a curve. Since 
$Y$ is regular it follows 
that $g(Y)$ is regular and hence a 
rational curve.
\end{proof}

\begin{remark} \em
In her unpublished thesis Renza Cortini 
has classified all  smooth surfaces $X$ of degree $d$ 
admitting a rational map $X\dasharrow \mathbb P^{2}$  of 
degree $d-2.$

\end{remark}

We can now improve \ref{stu}.

 \begin{proposition}\label{stu2}
 Let $X$ be a smooth surface in $\bP^{3}$ of degree $d,$  
 which contains no rational curves.
 Let $Y$ be a non-rational
  smooth surface with $p_g(Y)=q(Y)=0,$ and
let $ f: X \dasharrow Y$ be a generically 
  finite dominant rational  map of degree $m.$
  Assume moreover that $r=1,$
  in the notation of \S \ref{hurwitz}.
Then we have
$ K_{Y}^{2}=1$ 
 and $d\leq 6.$
\end{proposition}
\begin{proof}
 By \ref{cortini}, 1) we have $m\geq d-1$.
 Since $r=1$, 
from the proof of \ref{stu} we have
$$
 (d-1)K_{Y}^{2}\leq mK_{Y}^{2}=d+W^{2}\leq d.
 $$
This forces $K_{Y}^{2}=1$ and $m \leq d$. 
Now assume that $d>6.$
By \ref{cortini}, 2) it follows that $m=d.$
Recall the maps 
$h:Z\to Y$ and $\phi: Z \to X$ 
in diagram (\ref{diagram}).
From \ref{stu} again we have  
$W=0,$ hence:
 $$h^{\ast}K_Y=L=\phi^{\ast}H.$$  
 In particular we obtain $h_{\ast}L=dK_Y.$
Since $d>6$ the general fiber of the map $f:X \dasharrow Y$
is contained in a line. Thus
we obtain a surface $S$  in the Grassmannian of lines in $\mathbb P^{3}$ and 
a rational map $k: Y \dasharrow S$ birational onto its image.

We are going to show that a general plane $\Pi$ in $\mathbb P^{3}$
can be chosen in such a way that a number of special conditions are satisfied.
Let $S^0\subset S$ be the Zariski open subset of $S$
consisting of lines $\ell$ such that
 the cardinality of  $\ell\cap X$  is exactly $d=m$,
 and define:
$$U_0= \{\Pi\in { \mathbb P^{3}}^{\vee}: 
\Pi \supset\ell \mbox{\ for some\ } \ell \in S^0 \}. $$
 Let $\Gamma\subset X$ be the set of points of indeterminacy of the map
 $f,$ and define:
 $$
 U_1 =\{\Pi\in { \mathbb P^{3}}^{\vee}: \Pi\cap \Gamma=\emptyset \}.
 $$
 Define moreover:
 $$
  U_2 =\{\Pi \in  {\mathbb P^{3}}^{\vee}: \Pi \cap X
  {\rm \ is \ a \ connected \ smooth\  curve}\},
 $$ 
 $$
 U_3   =\{\Pi \in  {\mathbb P^{3}}^{\vee}:  
 \exists \ \ell \in S^0 , \exists \ x\in X:   \ell\cap \Pi=\{x\} \}.
 $$ 
We claim that the $U_i$  
are dense Zariski constructible subsets of ${\mathbb P^{3}}^{\vee}$,  so 
their intersection is non-empty.
This requires some well known basic facts,
except possibly the density of $U_{0}.$

This point is easily seen by means of projective duality.
It is enough to prove that every plane must contain some line of $S$.
There is a dual surface $S'$ in the Grassmannian of lines 
of ${\mathbb P^{3}}^{\vee}$ and the dual assertion is
that every point of ${\mathbb P^{3}}^{\vee}$
belongs to some line of $S'$.
Otherwise $S'$ covers a surface $X'$ in ${\mathbb P^{3}}^{\vee}$,
and through any two points of $X'$ there is a line of $S'$, so
$X'$ is a plane and $S'$ is the variety of lines in the plane.
In this case the dual surface $S$ is the variety of lines through
some fixed point in $\mathbb P^{3}.$ 
Then $Y$ is a rational surface, and this is a contradiction.

It follows that we can take  $\Pi \in U_0\cap U_1 \cap U_2 \cap U_3 .$
Since $\Pi\in U_2$  we see that
  $C= \Pi \cap X$  is a smooth curve of genus  $\frac {(d-1)(d-2)}{2}.$
Since   $ \Pi\in U_1$  we get that the restriction map $f_C: C \to f(C)$ 
is everywhere defined, and since
moreover  $ \Pi\in U_3$ then  $f_C$  is birational onto its image. 
 On the other hand since $\Pi \in U_0$ then $C$
 contains  $d$ distinct points,  of some line $\ell,$ which
 collapse on $f(C).$ Let   $a$ be  arithmetic 
 genus of $f(C)$. We have
 $$a\geq \frac {(d-1)(d-2)}{2}+
 \frac {d(d-1)}{2}=(d-1)^2.$$ 
We remark in fact that $f(C)$ is contained in
 a surface and has a $d-$ple point:
 the above bound follows from the adjunction 
 formula applied on the blow-up of $Y.$
 On the other hand since 
 $f(C)=h_{\ast}L=dK_Y$ we have
 $$2a-2= dK_Y \cdot (d+1)K_Y=d^2+d.$$
Hence
 $$ 2(d-1)^2 \leq d^2+d+2.$$
 That is $d^2\leq 5d$ and $d\leq 5.$
\end{proof}

To outline the importance of the number 
$s$ we give the following
\begin{definition}
We call the number $s=d-4-r$ the
 birational index of the ramification
of $f$ \label{birindex}.
\end{definition}

\subsection{General surfaces}
For a general surface in $\mathbb P^{3}$ 
all the assumptions required for the results 
in the present section are indeed satisfied.
We start by collecting some well known facts:

\begin{theorem} \label{factsonsurfaces}
Assume that $X$ is a general surface 
of $\mathbb P^{3}$ of degree 
$d\geq 5.$ Then
\begin{enumerate} \label{generale}
\item[i)] (Noether-Lefschetz) The N\'eron-Severi 
group of $X$ 
is generated by the hyperplane 
section $H.$
\item[ii)] (Lefschetz) The Hodge substructure of 
$H^{2}(X)$ orthogonal to the 
hyperplane section is irreducible.
\item[iii)] (Xu) The surface $X$ does not 
contain any rational or elliptic
curve.
\item[iv)] The only birational
 automorphism of $X$ is the identity.
\end{enumerate}
\end{theorem}
\begin{proof}
(i) follows from (ii), and for (ii) see Voisin \cite{Voisin}, \S 3.2.3.
(iii) is proved in \cite{Xu}, and (iv) is well known.
\end{proof}

We apply this to obtain:

\begin{proposition} \label{hodge}
Let  $X\subset \mathbb P^{3}$ 
be a general surface of degree 
$d\geq 5$ and $Y$ be a surface
 of general type. Let $f:X\dasharrow Y$
be a dominant rational
 map of degree $m>1.$
Then:
\noindent
\begin{enumerate}
\item $p_{g}(Y)=0,$ 
\item $Y$ is simply connected.
\end{enumerate}
\end{proposition}
\begin{proof}
\ 1.   Consider the Hodge structure 
map $f^{\ast}:H^{2}(Y) \to H^{2}(X)$,
which is defined by means of diagram (\ref{diagram})
as the composition of the ordinary pullback 
$h^*:H^{2}(Y) \to H^{2}(Z)$
followed by the Gysin map $\phi_*:H^{2}(Z) \to H^{2}(X)$,
and consider the injection 
$f^{\ast}:H^{2,0}(Y) \to H^{2,0}(X).$
Let $T_{Y}\supset H^{2,0}(Y) $ 
and $T_{X}\supset H^{2,0}(X) $  
be the  Hodge substructures orthogonal to the 
N\'eron-Severi Hodge structures
of $Y$ and respectively of $X.$
 We have $H^{2,0}(Y) =T_{Y}^{2,0}$ 
and $H^{2,0}(X) =T_{X}^{2,0}.$  
Then $f^{\ast}T_{Y}\subset T_{X}$ is a 
Hodge substructure of $ T_{X}.$ Assume, 
by contradiction, that $H^{2,0}(Y)\neq 0.$ Then
 $f^{\ast} H^{2,0}(Y)$ is not 
trivial and  hence by  \ref{generale} ii)  
the inclusion
$f^{\ast}T_{Y}\subset T_{X}$ is 
an equality. In particular 
$$f^{\ast}:H^{2,0}(Y) \to H^{2,0}(X)$$ is 
an isomorphism. It follows that
 the canonical map of $X$ factors
  through $f$ (as a rational map).
When $d>4$ the canonical map 
is an embedding. It would follow that $f$ 
is a birational map, which is a contradiction. 
\vskip 5pt

2.  Let 
 $\rho:W\to Y$ be the universal covering of $Y$.
Consider diagram (\ref{diagram}), in which $f$
is extended to a morphism $h$ on the blow up $Z$.
Since $Z$ is simply connected we may lift $h$ 
to a holomorphic map $g:Z\to W.$ It follows that
$g$ is surjective, $W$ is projective 
and the fundamental group
of $Y$ is finite.
 Since the deck 
 transformations of $\rho :W\to Y$
  give automorphisms of $W,$ and $X$ has only the
 trivial birational automorphism  it follows $\deg(g)>1.$
This defines a dominant rational map $g':X\dasharrow W$ 
of degree $>1.$
Using the first part of the proposition
we get  $p_g(W)=0.$ 
Since  $q(Y)=q(W)=0$
 we also obtain:
$$\chi (\cO_Y)= \chi (\cO_W)=1.$$ 
The proportionality theorem 
for the holomorphic
  Euler characteristic 
gives $$\chi (\cO_W)=\deg(\rho) \chi (\cO_Y).$$
Therefore  $\deg(\rho)=1,$ and the
 fundamental group of $Y$ is trivial.
\end{proof}

\begin{remark} \em A similar
 proposition holds for 
the general hypersurface 
of  $\mathbb P^{n}$
of degree $>n+1$.
\end{remark}

\section { Families of rational maps}

We consider a family of rational maps
from surfaces in $\mathbb P^{3}$ to surfaces of general type,
and we study the number of moduli  of the ramification divisors.
In the final subsection we prove our main result, 
that for small degree the general surface in $\mathbb P^{3}$ has no such map.

\subsection{Families and moduli} \label{famemod}
From now on $X$ will be a general 
smooth surface of $\mathbb P^3$ 
of degree $d\geq 5,$
 $Y$ will be a  simply
  connected  minimal surface of general type 
with $p_{g}(Y)=0.$
We also assume that $f:X \dasharrow Y$ is 
a dominant rational map of degree $m.$ 
The resolution of indeterminacy of $f$
is given in diagram (\ref{diagram}).
Moreover we would like to
 consider a family 
of such mappings.

Let $U$ be a smooth  base 
 variety  and $p\in U$ be
a general point. Assume we are given
a smooth family 
$q_1: \mathcal X \to U$
with $\cX \subset U \times \bP^{3},$
and smooth projective families
$q_2: \mathcal Z \to U$ and
$q_3: \mathcal Y \to U,$ 
and  a diagram of families
$$\begin{array}{ccc}
&\cZ& \\
&^{\scriptsize\Phi}\!\! \swarrow 
\hspace{10pt} \searrow \!\! ^{\scriptsize  H} \\
\cX \hspace{-10pt} & \underset{F}{\dasharrow} & \hspace{-10pt} {\cY}
\end{array}$$
 such that 
 $\Phi $ is birational 
 and ${\it H}$ is  generically finite and
 dominant.
 Moreover at the given point $p$ 
 the diagram of families specializes to 
 diagram (\ref{diagram}) relative to the map $f$.
 Thus we have defined a family of rational
  maps $F_t=  {\it H}_t \circ \Phi _t^{-1}.$
 
 We have now the moduli  
 space $\mathcal M_d$  of smooth surfaces 
 of degree $d$ in $\mathbb P^3$, and also  
 the connected component $\mathcal N$ of the
moduli space of minimal surfaces of general type
 with invariants $\chi(Y)$ and  $K^2_Y$,
 which contains all points $[Y_t=q_3^{-1}(t)]$ for $ t\in U.$   
  Define 
 $$N=\dim \mathcal N.$$
Consider the modular maps:
\begin{equation*} \label{modmap}
\mu: U \to \mathcal M_d \ \ \mbox{ and } \ \  \nu: U \to \mathcal N.
\end{equation*}
We assume that $\mu$
 is generically  finite and dominant, 
$$\dim U=M(d)=\binom{d+3}{3}-16.$$

For every rational map in the family there are several ramification loci.
Let $R_{t} \subset Z_{t}$  be the ordinary ramification divisor of $H_{t}$, 
and denote by $\overline{R}_{t}$ the support of $R_{t}$. 
Define
$B_{t} \subset Y_{t}$ to be the reduced branch divisor,
the support of the divisor ${H_{t}}_{*}(R_{t})$,
and define moreover $D_{t} \subset X_{t}$ to be the 
support of the divisor ${\Phi_{t}}_{*}(R_{t}),$
the part of the ramification that appears in $X_{t}$.
These divisors form
algebraic families, and those which are reduced divisors 
form algebraic families of curves in the sense of  \ref{modulardimension}
(shrinking the base U if necessary).
Let $\cR \subset \cZ$ be the relative ramification divisor of $H$ over $U,$ 
and denote by $\overline \cR \subset \cZ$
the total space of the family of curves $\overline R_{t}.$
Let $\cB \subset \cY$ be the total space of the family $B_{t},$ 
and let $\cD \subset \cX$ be the total space of the family $D_{t}.$ 
Here $\cB$ and $\cD$ are reduced.
Since $X_{t}$ is a general surface of $\bP^{3}$, then 
$D_{t}$ is a complete intersection curve, by \ref{factsonsurfaces} i). 
Moreover generically any component of $D_t$
has geometric genus $>1,$ by \ref{factsonsurfaces} iii).

We can write formulas
(\ref{can}) and (\ref{ram})
\noindent \begin{equation*} \label{cant}
 H_t^{\ast}(K_{Y_t})=rL_t-W_t,
 \end{equation*}
\begin{equation*} R_t= sL_t+ W_t +E_t. \label{ramt}
 \end{equation*}
 Note that the birational index $s$
  is  invariant in the family.
 We have $s+r=d-4$. It follows that $D_t$ 
 is a complete intersection of type 
 $(d,s')$ where $s'\leq s.$

We are going to study the modular dimension of this family $D_{t}$,
as defined in \ref{modulardimension}. 
An upper bound follows from considering this as a
family of complete intersection curves in $\bP^{3}.$

Define for $d > k$ and $d>2$ the function:
$$
M(d,k) = 
\begin{cases} \label{modcomint}
 \binom{d+3}{3} + \binom{k+3}{3}  - \binom{d-k+3}{3} -17  & \text{if } k>1\\
\binom{d+2}{2}-9& \text{if } k=1
\\ -1 & \text{if } k=0
\end{cases}.
$$
\begin{proposition} \label{complete}
Using the notation established above, we have:
$$M ({\mathcal D}/U) \leq M(d,s).$$
\end{proposition}
\begin{proof}
Let $\cD'$ be an irreducible component of $\cD$.
After a base change   $W \rightarrow U$
(and replacing $U$ by an appropriate open subset if necessary) 
 we have a family $\cC \to W$ whose members are
the normalizations $C_z$
of  irreducible components $A_{z}$  of  curves of the family ${\mathcal D}$, 
see \S \ref{moddim}. 
The curves in this family are complete intersections
of type $(d,s')$ with $s' \leq s$,
of geometric genus $g$. 
By \ref{factsonsurfaces} iii) we have that  $g>1$.
The case $s'=0$ is trivial and  $s'=1$
is similar to the case $s'>1$. So assume $s'>1$.

Let $\mathcal H(d,s')$ be the Hilbert scheme 
of complete intersection curves in $\bP^{3}$ of type $(d,s')$.
By acting with the projective 
group $G$ we find a family  of curves $hA_{z},$ $h\in G,$
of $\bP^3$  and this defines a morphism
$W \times G \to \mathcal H(d,s')$,
let $\widetilde W$ be the image of this morphism,
and consider the modular map
$\mu:\widetilde W \rightarrow {\mathfrak M}_{g}$.
Since $g>1$, the automorphisms group of $C_{z}$ is finite and hence
no continuous family of automorphisms of $\mathbb P^{3}$ can fix
the curve $A_{z}$.
It follows that  
$\dim \mu (\widetilde W) \leq \dim \widetilde W -\dim G=\dim \widetilde W -15$.
Therefore we have: 
$$M({\mathcal D'}/U) = \dim \mu (\widetilde W) \leq
\dim\widetilde W -15\leq
\dim {\mathcal H}(d,s') -15= M(d,s'),$$
and clearly $M(d,s') \leq M(d,s).$
\end{proof}

On the other hand,
a lower bound for the modular dimension comes from
the bend and break lemma.

\begin{proposition} \label{dimensioni}
Using the notation established above, we have:
$$M(d)-N \leq M(\cD/U)+1.$$
\end{proposition}

\begin{proof}
Let $\Gamma\subset U$ be the general fiber  of the modular map
$\nu:  U \to \mathcal N.$ 
We may assume that $\Gamma$ contains $p.$
Define $f = \dim \Gamma,$ and we remark that 
$$f \geq M(d)-N.$$
Let $\cX_{\Gamma}=q^{-1}_{1}(\Gamma) \ \ \mbox{and} \ \
 \cZ_{\Gamma}=q^{-1}_{2}(\Gamma)$  be the 
 restrictions of the families $\cX$ and $\cZ$ to $\Gamma$.
The maps parametrized by $\Gamma$ are
rational maps  $ F_t: X_t \dasharrow Y$ where
$Y$ is fixed. 
Let moreover $\cR_{\Gamma}$ and $\overline{\cR}_{\Gamma}$,
and $\cB_{\Gamma}$, $\cD_{\Gamma}$
be the restricted families of ramification
and branch divisors.

We recall from covering theory that covering maps of $Y$
with given branch locus $B$ and given degree $m$ 
are classified in terms of homomorphisms
from the fundamental group of $ Y\setminus B$
to the symmetric group in $m$ elements 
(see \cite{miranda} ch. III \S 4).
It follows that the map $t \mapsto B_{t}$ from $\Gamma$
to the Hilbert scheme of $Y$ is generically finite
onto its image (since so is the map $U \to \mathcal M_{d}$
in our assumptions).
So the branch loci $B_t$, $t\in\Gamma,$
form a family of curves in $Y$ of dimension $f.$ 
It follows from \ref{curve} that the 
modular dimension of this family is
$M(\cB_{\Gamma}/\Gamma)\geq f-1$.
A fortiori the family of 
 ramification divisors $R_t$, $t\in\Gamma,$
has modular  dimension 
$M(\overline{\cR}_{\Gamma}/\Gamma) \geq M(\cB_{\Gamma}/\Gamma)\geq f-1.$

Assume now that $s>0$.
The family $D_{t},$ $t\in\Gamma,$
has the same modular dimension
$M(\cD_{\Gamma}/\Gamma) = M(\overline{\cR}_{\Gamma}/\Gamma).$
An irreducible component of $\overline{\cR}_{\Gamma}$
whose fibres consist of rational curves
has modular dimension 0.
Any other irreducible component
arises from an irreducible component of $\cD_{\Gamma}$
with the same modular map.
Summing up we have:
$$f-1 \leq M(\cB_{\Gamma}/\Gamma) 
\leq M(\overline{\cR}_{\Gamma}/\Gamma) 
= M(\cD_{\Gamma}/\Gamma) \leq  M(\cD/U).$$ 

Finally, if $s=0$ then $D_{t}=0$ and $B_{t}$ consists of rational curves,
which  do not move in $Y$, 
hence it follows that $f=0,$ and this gives the statement.
\end{proof}

\subsection{Main result}

We are in a position  to prove our main theorem,
by combining the previous results.
We keep the notation of the last section.

From \ref{complete} and  \ref{dimensioni} 
we obtain the inequality: 
\begin{equation} \label{inequality}
M(d)-N \leq M(d,s)+1,
\end{equation}
which is a necessary condition for
the existence of a rational map $f,$
with birationality index $s,$ 
on a general surface $X$ of degree $d$ in $\bP^{3}.$
Moreover from \ref{mo2} we have:
$$N \leq 19.$$

\begin{theorem} \label{main} If  $5 \leq d\leq 11$ then the general 
surface of degree $d$ has no (non-trivial)
rational map which dominates a surface of general type.
\end{theorem}
\begin{proof}
If $s=0$ then inequality (\ref{inequality})
gives  $M(d)-N \leq 0$, that is: 
$$\binom {d+3}{3}-16\leq N;$$ since $d>4$
 this gives $N\geq 40$, a
  contradiction.
If $s=1$ since $r>0$
 we have $d\geq 6$. From  (\ref{inequality})
 we have the inequality:
$$\binom {d+3}{3}-16-\binom {d+2}{2}+9\leq N+1;$$
 for $d>5$ this gives $N\geq 48$, 
 again a contradiction. 
Assume $s>1$, so that  $d>6$. By \ref{stu2} 
we have $r>1.$
From  (\ref{inequality}) we have:
$$\binom {d+3}{3}-16-N\leq  \binom{d+3}{3} +
 \binom{s+3}{3}  - \binom{d-s+3}{3} -17+1 $$
whence we obtain: 
$$19 \geq  N \geq \binom{d-s+3}{3} - 
\binom{s+3}{3}=\binom{r+7}{3} - 
\binom{s+3}{3}.$$
Now since $r\geq 2$ 
we have:
$$19\geq  N\geq \binom{9}{3} - 
\binom{s+3}{3}=84-\binom{s+3}{3},$$
and so: 
$$ \binom{s+3}{3}\geq  65.$$
This gives $s\geq 6$ and hence $d=s+r+4 \geq 12.$ 
\end{proof}

\begin{remark} \em
The previous computation 
proves that the only possible case
for $d=12$ gives   $s=6$  and $r=2.$  
\end{remark}

\vskip 40pt
\noindent {\it Acknowledgements:} 
It is a pleasure to thank Shigeru Iitaka and Fabrizio Catanese
for their interest in the present work. We are also grateful to
Alessandro Ghigi, Roberto Pignatelli
and Enrico Schlesinger for many fruitful conversations.
Finally we thank the referees for several useful suggestions and remarks.

\bigskip

\noindent {\sc Lucio Guerra }\\
Dipartimento di Matematica, Universit\`a di Perugia\\
Via Vanvitelli 1, 06123  Perugia, Italia\\
{\tt guerra@unipg.it}
\bigskip

\noindent {\sc Gian Pietro Pirola}\\
Dipartimento di Matematica, Universit\`a di Pavia\\
via Ferrata 1, 27100 Pavia, Italia\\
{\tt gianpietro.pirola@unipv.it}


\vskip 40pt


\begin{thebibliography}{11}

 
  \bibitem{Amerik}
  E.~Amerik.
 \newblock {\em {O}n a problem of {N}oether-{L}efschetz type.}
 \newblock { Compositio Math. {\bf 112} 
(1998), no. 3, 255--271.}

\bibitem{Barlow}
R. ~Barlow.
{\em A simply connected surface of general type with $p\sb g=0$.}
Invent. Math. {\bf 79} (1985), no. 2, 293--301. 

\bibitem{BHPV}
W. ~Barth, K. ~Hulek, C. ~Peters, A. ~Van de ~Ven. 
{\em Compact complex surfaces.}
Second edition. Ergebnisse der Mathematik und ihrer Grenzgebiete. 3. Folge, 4.
Springer-Verlag, Berlin, 2004. 

\bibitem{bogomolov}
F.~Bogomolov. 
 \newblock {\em Holomorphic tensor and vectors bundles on projetive varieties.} 
 \newblock {Math. USSR Izv. {\bf 13} (1979, 499--555.}

\bibitem{Catanese}
F.~Catanese. 
 \newblock {\em On the moduli spaces of surfaces of general type.} 
 \newblock {J. Diff. Geom. {\bf 19} (1984), no. 2, 483--515.}
 
\bibitem{Catanese2}
F. ~Catanese. 
{\em Moduli of algebraic surfaces.}
Theory of moduli (Montecatini Terme, 1985), 1--83,
Lecture Notes in Math., 1337,
Springer, Berlin, 1988. 

\bibitem{Cheltsov}
I.~Cheltsov. 
 \newblock {\em Points in projective spaces and applications.} 
 \newblock {  math.AG/05115778.}
 
 \bibitem{ciro}
C. ~Ciliberto, G. van der ~Geer.
{\em On the Jacobian of a hyperplane section of a surface.} 
Classification of irregular varieties (Trento, 1990), 
Lecture Notes in Math., 1515, 
Springer, Berlin, 1992, 33--40.
  
  \bibitem{enoki}
  I.~Enoki.
 \newblock {\em {S}tability and negativity for tangent shaeaves
 of minimal K\"{a}hler spaces.}
 \newblock {Geometry and analysis on manifolds (Katata/Kyoto, 1987), 118--126,
Lecture Notes in Math., 1339,
Springer, Berlin, 1988.}

\bibitem{freedman} M.~Freedman.
{\em The topology of 4-manifolds.} J. Diff. Geom. {\bf 17} (1982), 357-454.

 
\bibitem{griffithsha}
P.~Griffiths and J.~Harris.  
\newblock {\em Principles of {A}lgebraic {G}eometry}.
\newblock {Wiley-Interscience, New York,  { 1978}.}
 
 \bibitem{Guerra}
L.~Guerra. 
 \newblock {\em Complexity of Chow varieties and number of 
 morphisms on surfaces of general type.} 
 \newblock {Manuscripta Math. {\bf 98}  (1999), no. 1, 1--8.}
 
 \bibitem{lggpp}
L.~Guerra and G.~Pirola.
 \newblock {\em  On the finiteness theorem for rational maps on a variety of general type.} 
 \newblock {in preparation.}
  

\bibitem{Haconk}
C.~Hacon and J.~McKernan.
\newblock {\em Boundedness of pluricanonical maps
of varieties of general type}.
 \newblock{Invent. Math.  {\bf 166} (2006), no. 1, 1--25.}
 
 \bibitem{Heier}
G.~Heier. 
\newblock{\em Effective finiteness theorems for maps between
canonically polarized compact complex manifolds}.
 \newblock{Math. Nachr. {\bf 278} (2005), no. 1-2, 133--140.}
 
 \bibitem{Kani}
E.~Kani.
 \newblock {\em Bounds on the number of nonrational subfields
of a function field}. Invent. Math  \newblock{ {\bf 85} (1986), no. 1, 185--198}.
  
 \bibitem{KO}
S. ~Kobayashi, T. ~Ochiai.
{\em Meromorphic mappings onto compact complex spaces of general type}.
Invent. Math. {\bf 31} (1975), no. 1, 7--16.

 \bibitem{Kobayashi}
S. ~Kobayashi . 
\newblock {\em Differential Geometry  of Complex Vector Bundles}
\newblock{ Publ. of Math Soc. of Japan and Princton University Press}.
 (1987).


\bibitem{LeePark}
Y.~ Lee and J. ~Park.
{\em A simply connected surface of general type with $p_g=0$ and $K^2=2$.}
 ArXiv: math.AG/0609072.

  \bibitem{lopezpirola}  
A.~Lopez and G.~Pirola.
 \newblock {\em On the curves through a general point of a smooth
  surface in $\mathbb P\sp 3$.} 
 \newblock { Math. Z. {\bf 219} (1995), no. 1, 93--106.}

  \bibitem{maehara}
K.~Maehara. 
 \newblock {\em A finiteness property of varieties of general type.} 
 \newblock {Math. Ann. {\bf 262} (1983), no. 1, 101--123.}
 
 \bibitem{miranda}
R.~Miranda.
\newblock{\em Algebraic curves and Riemann surfaces. }
Graduate Studies in Mathematics, 5. 
American Mathematical Society, Providence, RI, 1995. 

\bibitem{Mori}
S.~Mori.
{\em Projective manifolds with ample tangent bundles}.
Ann. of Math. (2) {\bf 110} (1979), no. 3, 593--606.

\bibitem{GIT}
D.~Mumford, J.~Fogarty.
\newblock{\em Geometric invariant theory.}
Second edition. Ergebnisse der Mathematik und ihrer Grenzgebiete, 34. 
Springer-Verlag, Berlin, 1982. 

\bibitem{NaranjoPirola}
J.~Naranjo and G.~Pirola.  \newblock {\em Bounds of the number of
rational maps betweeen varieties of general type}. 
 \newblock{ArXiv: math.AG/0511463}.
 
\bibitem{reid}
M.~Reid. 
{\em Bogomolov's theorem $c_1^2\leq 4c_2$.}  
Proc. Symp. Algebraic Geometry (Kyoto, 1977), 623-642.

 \bibitem{Severi}  F.~Severi.
 \newblock {\em  Le corrispondenze fra i punti di una curva variabile
  sopra una superficie algebrica},
Math Ann. {\bf 74} (1913), 511-544.


  \bibitem{sugiyama}
 K.~Sugiyama.
 \newblock {\em {O}n tangent sheaves of minimal varieties.}
 \newblock {K\"{a}hler metrics and Moduli spaces,  85--103, 
 Adv. Stud.  Pure Math., 18-{\rm II},
Academic Press, Boston, MA, 1990.}

\bibitem{takayama}
S.~Takayama. 
{\em Pluricanonical systems on
algebraic varieties of general type.}  Invent. Math. {\bf 165} (2006), no. 3, 551--587.
 
\bibitem{Tanabe}
M.~Tanabe.  \newblock {\em Bounds on the number of holomorphic maps of compact Riemann surfaces}.
 \newblock{Proc. Amer. Math. Soc.  
 \newblock{{\bf133} (2005), no. 10, 3057--3064.}}

\bibitem{Tsuji}
H.~Tsuji. 
\newblock{\em Pluricanonical systems of projective varieties of general type II}.
 \newblock{ ArXiv: math.CV/0409318.}
 
\bibitem{Voisin}
C. Voisin.
{\em Hodge theory and complex algebraic geometry}. II.
Cambridge Studies in Advanced Mathematics, 77.
Cambridge University Press, Cambridge, 2003. 

\bibitem{yau}
S.~Yau. 
{\em Calabi's conjecture and some new results in algebraic geometry.}  Proc.. Nat. Ac. Sc. Usa . {\bf 74} (1977), 1798-1799.

\bibitem{yau2}
S.~Yau. 
{\em On the Ricci curvature of a complete Kaehler manifold and the complex Monge-Amp\'ere equation.} Comm. Pure Appl. Math. {\bf 31} (1978), 339-411.
 
  \bibitem{Xu}  
G.~Xu.
\newblock {\em Subvarieties of general hypersurfaces in projective space}. 
\newblock {J. Differential geom. {\bf 39} (1994), no. 1, 139-179.}

 
\end{thebibliography}
\end{document}